\documentclass[10pt]{amsart}
\usepackage{amsmath}
\usepackage{amsfonts}
\usepackage{amssymb}

\newcommand{\Ext}{\operatorname{Ext}}
\newcommand{\End}{\operatorname{End}}
\newcommand{\Hom}{\operatorname{Hom}}
\newcommand{\F}{{\mathcal F}}
\newcommand{\T}{{\mathcal T}}
\newcommand{\tT}{\tilde{{\mathcal T}}}
\newcommand{\Y}{{\mathcal Y}}
\newcommand{\G}{{\mathcal G}}
\newcommand{\ks}{k\Sigma_d}
\newcommand{\soc}{\operatorname{soc}}
\newcommand{\mo}{\operatorname{mod-}}
\newcommand{\di}{\operatorname{dim}}
\newcommand{\sgn}{\operatorname{sgn}}
\newcommand{\rad}{\operatorname{rad}}
\newcommand{\ind}{\operatorname{Ind}}
\newcommand{\cha}{\operatorname{char}}

\numberwithin{equation}{section}

\theoremstyle{plain}
\newtheorem{theorem}{Theorem}[section]
\newtheorem{cor}[theorem]{Corollary}
\newtheorem{lemma}[theorem]{Lemma}
\newtheorem{prop}[theorem]{Proposition}

\newtheorem{conj}[theorem]{Conjecture}

\newtheorem{example}[theorem]{Example}
\newtheorem{remark}[theorem]{Remark}
\newtheorem{problem}[theorem]{Problem}

\begin{document}
\title[Specht filtrations]
{\bf Symmetric group modules with Specht and dual Specht filtrations}
\author{\sc David J. Hemmer}

\address
{University of Toledo\\ Department of Mathematics \\
2801 W. Bancroft\\Toledo\\ OH~43606, USA}
\thanks{Research of the author was supported in part by an NSA
Young Investigator's Grant } \email{david.hemmer@utoledo.edu}

\date{April 2006}

\begin{abstract} The author and Nakano recently proved that multiplicities in a Specht filtration of a symmetric group module are well-defined precisely when the characteristic is at least five.  This result suggested the possibility of a symmetric group theory analogous to that of good filtrations and tilting modules for $GL_n(k)$. This paper is an initial attempt at such a theory. We obtain two sufficient conditions that ensure a module has a Specht filtration, and a formula for the filtration multiplicities. We then study the categories of modules that satisfy the conditions, in the process obtaining a new result on Specht module cohomology.

Next we consider symmetric group modules that have both Specht and dual Specht filtrations. Unlike tilting modules for $GL_n(k)$, these modules need not be self-dual, and there is no nice tensor product theorem. We prove a correspondence between indecomposable self-dual modules with Specht filtrations and a collection of $GL_n(k)$-modules which behave like tilting modules under the tilting functor. We give some evidence that indecomposable self-dual symmetric group modules with Specht filtrations may be self-dual trivial source modules.

\end{abstract}

\maketitle

\section{Introduction} Let $k$ be an algebraically closed field and let $G$ be a reductive algebraic group over $k$. A rational $G$-module is said to have a good filtration if it has a filtration with successive quotients
isomorphic to induced modules $\nabla(\lambda)$. There is a simple cohomological criterion for having a good filtration, and a formula for the multiplicities, which are
independent of the choice of filtration. The indecomposable modules which have both a good and a Weyl filtration are called tilting modules, and are labelled by
dominant weights.

Until recently it was thought no such theory could exist for Specht and dual Specht filtrations of symmetric group modules, since well-known examples in characteristic two and three demonstrated that filtration multiplicities are not even well-defined. Nevertheless, in  \cite{HNspechtfilt} it was shown that the multiplicities are well-defined as long as $\cha k>3$. Necessary and sufficient conditions for determining if a module has a Specht and/or dual Specht filtration were obtained. However the conditions are not in terms of symmetric group cohomology, but rather are stated in terms of $GL_n$ cohomology and the adjoint Schur functor.

This paper is a first attempt at a theory of Specht filtrations. We give two different sufficient conditions for a $k\Sigma_d$ module to have a Specht (or dual Specht) filtration. Although the conditions are not necessary, they have the advantage
of being stated entirely in terms of the symmetric group theory. For modules satisfying the conditions, we obtain a formula for the filtration multiplicities which generalizes a known formula for Young modules.  The collection of modules satisfying the condition gives an interesting subcategory of mod-$k\Sigma_d$, which we study.

We then  consider modules which have both Specht and dual-Specht filtrations. We demonstrate that they are not as well-behaved as tilting modules. Unlike tilting modules, they need not be self-dual. The tensor product of two modules which have both Specht and dual Specht filtrations may have neither! We believe however that a classification of the indecomposable {\it self-dual} modules with both filtrations is possible. We show they are in correspondence with $GL_n(k)$ modules satisfying a certain natural property under the tilting functor. We give some evidence, inspired by recent work of Paget and Wildon, to suggest these modules may be the self-dual trivial-source modules.

\section{Notation and Preliminaries}

We assume $\cha k>3$ throughout, and emphasize that many of the results do not otherwise hold. Let $G=GL_n(k)$ and $V\cong k^n$ be the natural $G$-module.  Let $$S:=S(n,d)\cong \End_{\ks}(V^{\otimes d})$$ be the Schur algebra. Modules for $S$ correspond to polynomial $G$ modules of homogeneous degree $d$. Our basic references for representation theory of $k\Sigma_d$ and of $S(n,d)$ are \cite{Jamesbook} and \cite{Martinbook}. 

Write $\lambda \vdash d$ for $\lambda=(\lambda_1, \ldots , \lambda_t)$ a partition of $d$ and let $\Lambda^+(d)$ be all partitions of $d$. Let $\Lambda^+(n,d)$ denote the partitions of $d$ with at most $n$ parts. Let $\Lambda(n,d)$ denote the compositions of $d$ with at most $n$ parts and let $\lambda'$ denote the transpose of $\lambda$. Simple $S$-modules are indexed by $\Lambda^+(n,d)$ and denoted $L(\lambda)$. Let $\nabla(\lambda)$ and $\Delta(\lambda)$ denote the induced and Weyl modules, and $P(\lambda)$, $I(\lambda)$ and $T(\lambda)$ the projective, injective  and tilting modules for $S$. The tensor products of symmetric (resp. exterior) powers of $V$ are denoted $S^\lambda(V)$ (resp. $\Lambda^\lambda(V))$. Descriptions of these modules can be found in \cite{Martinbook}. Let $\F(\nabla)$  (resp. $\F(\Delta)$) be the set of $S$-modules having good (resp. Weyl) filtrations.

A partition $\lambda$ is called $p$-restricted if $\lambda_i - \lambda_{i+1}<p$ for all $i$. It is $p$-regular if $\lambda'$ is $p$-restricted. The simple $\ks$ modules are indexed by $p$-restricted partitions and denoted $D_\lambda$. They can also be indexed by $p$-regular partitions and denoted $D^\mu$.  The Specht, Young and permutation modules for $\ks$ are indexed by $\Lambda^+(d)$ and denoted $S^\lambda$, $Y^\lambda$ and $M^\lambda$. Let $S_\lambda = (S^\lambda)^*$. Recall that:

\begin{eqnarray}S^\lambda \otimes \sgn \cong S_{\lambda '}&&\\\nonumber
D^\lambda \otimes \sgn \cong D_{\lambda'}.&&
\label{tensorsignwithspecht}
\end{eqnarray}

\subsection{Schur and adjoint Schur functors}

For $n \geq d$ let $e\in S(n,d)$ denote the idempotent described in \cite[(6.1)]{Greenpolygln}. Then $eSe \cong k\Sigma_d$, and the Schur functor $\F: \mbox{\rm mod-}S \rightarrow \mbox{\rm mod-}k\Sigma_d$ is defined by $\F(U):=eU$. Let $\tau$ denote the usual contravariant duality on $\mo S$. Then $\F$ is compatible with $\tau$ and
the usual duality on $\mo  \ks$, i.e. $\F(U^\tau) \cong (\F(U))^*$.

The Schur functor is an exact, covariant functor with:

\begin{equation}
\label{imageofSchur}
\begin{array}{lll}
  \mathcal{F}(\nabla(\lambda))= S^\lambda &
  \mathcal{F}(\Delta(\lambda)) = S_\lambda &
  \mathcal{F}(L(\lambda)) = D_\lambda \mbox{\rm \, or 0}\\[0.1in]
  \mathcal{F}(P(\lambda))=Y^\lambda &
  \mathcal{F}(I(\lambda))=Y^\lambda &
    \mathcal{F}(T(\lambda))=Y^{\lambda'} \otimes \sgn.\\
\end{array}
\end{equation}

\noindent The Schur functor admits a right adjoint functor $\G:
\mbox{\rm mod-}k\Sigma_d \rightarrow \mbox{\rm mod-}S$ defined by:
\begin{eqnarray*}\G(N)&:=&\Hom_{k\Sigma_d}(V^{\otimes
d},N)\\
&\cong & \Hom_{eSe}(eS,N). \end{eqnarray*} The functor $\G$ is a one-sided inverse to $\F$, i.e. $\F (\G(M)) \cong M$. The functor $\G$ is
only left exact, and so has higher right derived functors:
$$R^i\G(N) = \Ext^i_{k\Sigma_d}(V^{\otimes d},N).$$

We now collect some known results about $\G$ and $R^1\G$. In \cite{DENspectral} a Grothendieck spectral sequence is constructed to relate the cohomology of $S$ and $k\Sigma_d$ using $\F$ and $\G$. We will only use the related five-term exact sequence which begins:

\begin{equation}
\label{fiveterm} 0 \rightarrow \Ext^1_S(U,\G(N)) \rightarrow
\Ext^1_{\ks}(\F(U),N) \rightarrow \Hom_S(U,R^1\G(N))\rightarrow \cdots .
\end{equation}

\noindent Recall that we are assuming $p>3$ throughout, indeed most of the results below fail for $p\leq 3$.

\begin{prop}
\label{propertiesofG}
\begin{itemize}
\item[]
 \item[(i)]\cite[3.2]{KNcomparingcohom}.
$\G(S_\lambda) \cong \Delta(\lambda).$

\item[(ii)]\cite[5.2.4]{CPSStratmemoirs} $\G(Y^\lambda \otimes
\sgn) \cong T(\lambda')$, $\quad \G(M^\lambda \otimes \sgn) \cong \Lambda^\lambda(V).$

\item[(iii)]\cite[3.8.2]{HNspechtfilt} $\G(Y^\lambda) \cong
P(\lambda), \quad \G(M^\lambda) \cong S^\lambda(V).$

\item[(iv)]\cite[3.5.1]{HNspechtfilt} $M \in \mo \ks$ has a dual
Specht filtration if and only if $\G(M)$ has a Weyl filtration.

\item[(v)]\cite[2.1.2]{DENspectral} $\Hom_S(L(\lambda),\G(M))=0$
unless $\lambda$ is $p$-restricted. In particular if $L(\mu)$ lies in ${\rm
soc}(\Delta(\lambda))$ then $\mu$ must be $p$-restricted.

\end{itemize}
\end{prop}
\noindent We will also need some information about $R^1\mathcal{G}$:

\begin{prop}
\label{R1Gvanishesondualspechtmodules}
\begin{itemize}
\item[] \item[(i)]\cite[6.4]{KNcomparingcohom} $R^1\G(S_\lambda)=0.$

\item[(ii)] If $M$ has dual Specht filtration then
$\Ext^1_S(U,\G(M)) \cong \Ext^1_{\ks}(\F(U),M)$.

\end{itemize}
\end{prop}

\begin{proof}
Notice that (ii) follows from (i) and (\ref{fiveterm}).

\end{proof}

\section{A filtration criterion and multiplicity formula}
\label{filtrationcriterionsection}
 There is a well-known necessary and sufficient condition for an $S$ module to have a good filtration and a formula for the filtration multiplicities.

\begin{prop} \cite[Prop. A2.2]{DonkinqSchurbook} Let $V \in \mo S(n,d)$. Then:
\label{goodweylfiltrationcriterion}
\begin{itemize}
\item[(i)]$V\in \F(\Delta)$ if and only if $\Ext^1_S(V, \nabla(\lambda))=0$ for all $\lambda \in \Lambda^+(n,d)$. If so then the multiplicity of $\Delta(\lambda)$ is independent of the choice of Weyl filtration and given by:
$$[V:\Delta(\lambda)]=\di_k\Hom_S(V,\nabla(\lambda)).$$

\item[(ii)]$V \in \F(\nabla)$ if and only if $\Ext^1_S(\Delta(\lambda),V)=0$ for all $\lambda \in \Lambda^+(n,d)$. If so then the multiplicity of $\nabla(\lambda)$ is independent of the choice of good filtration and given by:
$$[V:\nabla(\lambda)]=\di_k\Hom_S(\Delta(\lambda),V).$$

\end{itemize}
\end{prop}

Propositions \ref{propertiesofG}(iv) and \ref{goodweylfiltrationcriterion} immediately imply that a $k\Sigma_d$ module $M$  has a dual Specht filtration if and only if $\Ext^1_S(\G(M), \nabla(\lambda))=0$ for all $\lambda \in \Lambda^+(d)$. It seems unsatisfying to have a condition stated in terms of the $S$-cohomology of $\G(M)$ instead of the $\Sigma_d$-cohomology of $M$. The corresponding multiplicity formula is also in terms of $S$, namely $$[M:S_\mu]=\di_k\Hom_S(\G(M), \nabla(\mu)).$$
Given a $\ks$ module with a Specht filtration, it is natural to ask if the filtration multiplicities are given by the dimension of some space of $k\Sigma_d$-homomorphisms. There is one situation where this is known to be the case, namely for Young modules in characteristic $p>3$. It is well known that $Y^\lambda$ is self-dual with both a Specht and dual Specht filtration. The filtration multiplicities are well-defined in characteristic $p>3$ by the main result of \cite{HNspechtfilt}, and they are given by the following:

\begin{prop}
\label{propyoungfiltmultiplicites}Let $[Y^\lambda:S_\mu]$ denote the multiplicity of $S_\mu$ in a dual Specht filtration of $Y^\lambda$. Then:
\begin{align}
\label{eqdecequal}
[Y^\lambda:S_\mu]&=\di_k\Hom_{k\Sigma_d}(Y^\lambda, S_\mu)\\
\notag
&= \di_k\Hom_{k\Sigma_d}(S^\mu,Y^\lambda).
\end{align}
\end{prop}
\begin{proof}We show both sides of \eqref{eqdecequal} are given by the decomposition number $[\Delta(\mu):L(\lambda)]$. First observe that

\begin{align*}
[\Delta(\mu):L(\lambda)]&=\di_k\Hom_S(P(\lambda),\Delta(\mu))\\
&=\di_k\Hom_S(P(\lambda),\G(S_\mu))\\
&=\di_k\Hom_{k\Sigma_d}(Y^\lambda,S_\mu)\quad \text{by adjointness.}\\
\end{align*}

\noindent But  the decomposition number is also a filtration multiplicity in a tilting module:
\begin{align*}
[\Delta(\mu):L(\lambda)] &=[T(\lambda'):\nabla(\mu')]\quad \text{by \cite[Lemma 3.1]{Donkingtiltingalgebraic}}\\ 
&= [Y^\lambda \otimes \sgn: S^{\mu'}]\quad \text{since $\F$ is exact}\\
&=[Y^\lambda: S_\mu]\quad \text{by \eqref{tensorsignwithspecht}}. \qedhere\\
\end{align*} 
\end{proof}

More generally we can ask:
\begin{problem}
\label{problemaboutformulaformult}
Can one classify which $\ks$ modules with Specht or dual Specht filtrations have multiplicities given by a formula like that in Prop. \ref{propyoungfiltmultiplicites}?
\end{problem}

In Section \ref{subsectionsymmgroupcriterion} we give a large class of modules for which Problem \ref{problemaboutformulaformult} has an affirmative answer.

\subsection{A symmetric group filtration criterion}
\label{subsectionsymmgroupcriterion}

We  now prove that the  symmetric group conditions corresponding to Prop. \ref{goodweylfiltrationcriterion} under the Schur functor are sufficient to guarantee existence of filtrations. We obtain a pair of conditions which  guarantee a Specht filtration and a pair which guarantee a dual Specht filtration. These are the first known conditions stated in terms of symmetric group cohomology that guarantee modules have certain filtrations. We also obtain a formula for the filtration multiplicities which generalizes Proposition \ref{propyoungfiltmultiplicites}. Notice that (i) and (iv) below correspond to Prop.
\ref{goodweylfiltrationcriterion} under $\F$.
\begin{theorem}
\label{mainfilttheorem} Let $M \in \mo
\Sigma_d$. 
\begin{itemize}
\item[(i)]If $\Ext^1_{\ks}(M,S^\lambda)=0 \,\,\,\forall \lambda
\in \Lambda^+(d)$ then $M$ has a dual Specht filtration. The multiplicity of $S_\mu$ in any such filtration is given by $\di_k\Hom_{k\Sigma_d}(M,S^\mu)$

\item[(ii)]If $\Ext^1_{\ks}(S^\lambda,M)=0 \,\,\,\forall \lambda
\in \Lambda^+(d)$ then $M$ has a dual Specht filtration. The multiplicity of $S_\mu$ in any such filtration is given by $\di_k\Hom_{k\Sigma_d}(S^\mu,M)$

\item[(iii)]If $\Ext^1_{\ks}(M,S_\lambda)=0 \,\,\,\forall \lambda
\in \Lambda^+(d)$ then $M$ has a Specht filtration. The multiplicity of $S^\mu$ in any such filtration is given by $\di_k\Hom_{k\Sigma_d}(M,S_\mu)$

\item[(iv)]If $\Ext^1_{\ks}(S_\lambda,M)=0\,\,\, \forall \lambda
\in \Lambda^+(d)$ then $M$ has a Specht filtration. The multiplicity of $S^\mu$ in any such filtration is given by $\di_k\Hom_{k\Sigma_d}(S_\mu,M)$.
\end{itemize}
\end{theorem}
\begin{proof}
We first prove (iv), so assume $\Ext^1_{\ks}(S_\lambda,M)=0\,\,\,
\forall \lambda \in \Lambda^+(d)$. Then by (\ref{fiveterm}) we get
that $\Ext^1_S(\Delta(\lambda),\G(M))=0 \,\,\, \forall \lambda \in
\Lambda(d)$. By Prop. \ref{goodweylfiltrationcriterion},
$\G(M)$ has a good filtration, and the multiplicity of $\nabla(\mu)$ in $\G(M)$ is given by $\Hom(\Delta(\mu),\G(M))$.  So $M=\F(\G(M))$ has a  Specht
filtration and the multiplicity formula follows from the adjointness of $\F$ and $\G$.  Now (i) follows immediately. To obtain (ii) and (iii)
 use \eqref{tensorsignwithspecht} and Prop. 3.1.

\end{proof}

We remark that the two criteria in each pair detect different modules. For instance a Young module $Y^\lambda$ has both a Specht and dual Specht filtration. It is known that  $\Ext^1_{\ks}(S^\lambda, Y^\mu)$ is always zero but $\Ext^1_{\ks}(Y^\mu, S^\lambda)$ may be nonzero \cite[4.1]{HNspechtfilt}, so the dual Specht filtration of $Y^\mu$ is detected by Theorem   \ref{mainfilttheorem}(ii) but not (i). Note that the multiplicity formula in Theorem \ref{mainfilttheorem}(ii) is a significant generalization of Prop. \ref{propyoungfiltmultiplicites}.

The criteria in Theorem \ref{mainfilttheorem} are sufficient but not necessary; it is  possible for a module to have a Specht filtration but not satisfy either Theorem \ref{mainfilttheorem}(iii) or (iv). The structure of the Young and Specht modules for $\Sigma_{2p}$ is well understood. (see for example \cite{Martinprojfor2p}) From this one can easily calculate:

\begin{example}
\label{missed}
 Let $p=5$ and $d=10$. Then
$$\Ext^1_{\Sigma_{10}}(S^{531^2},S_{51^5}) \cong k \cong
\Ext^1_{\Sigma_{10}}(S_{55},S^{531^2})$$
\end{example}
\noindent so the  Specht filtration of $S^{531^2}$ is not detected by either Theorem \ref{mainfilttheorem}(iii) or (iv).

\subsection{A covariantly finite subcategory of $\mo \ks$}

The categories of $k\Sigma_d$ modules which satisfy the various conditions in Theorem \ref{mainfilttheorem} should be interesting  to study. Let $\Theta=\{S^\lambda \mid \lambda \in \Lambda^+(d)\}$ and let $\F(\Theta)$ denote the full subcategory of $\mod k\Sigma_d$ having Specht filtrations. $\F(\Theta)$ is obviously closed under extensions and it follows from \cite[3.6.2]{HNspechtfilt} that it
is closed under direct summands. It is also the case that:

\begin{lemma}
\label{spechtext} \cite[4.2.1]{HNspechtfilt} Let $p>3$ and suppose
$\mu \not \!\rhd \lambda$. Then $\Ext^1_{\ks}(S^\mu,S^\lambda)=0$.
\end{lemma}

\noindent Lemma \ref{spechtext} together with the main result of
\cite{Ringcatmodalmostsplit} imply that $\F(\Theta)$ is
functorially finite and has almost split sequences.

Observe that  duality interchanges modules satisfying
condition (i) with condition (iv) and (ii) with (iii) in Theorem \ref{mainfilttheorem}. Tensoring with $\sgn$ interchanges those satisfying (i) with (iii) and (ii)
with (iv). Thus if we want to study the module categories which arise from Theorem \ref{mainfilttheorem}, we may without loss of generality consider Theorem \ref{mainfilttheorem}(ii) and study the following category:
\[
{\mathcal Y}(\Theta):= \{M \in \mo k\Sigma_d \mid
\Ext^1_{\ks}(S^\lambda,M)=0\,\, \forall \,\, \lambda \in
\Lambda^+(d)\}.
\]
Notice that Theorem \ref{mainfilttheorem}(ii) guarantees all modules in $\Y(\Theta)$ have dual Specht filtrations. Categories such as $\Y(\Theta)$ are studied in \cite{Ringcatmodalmostsplit} where, for example, it is shown that they are covariantly finite, i.e. every $\ks$-module has a left ${\mathcal Y}(\Theta)$-approximation. These approximations are constructed explicitly in \cite[Sect. 2]{Ringcatmodalmostsplit}.
For another example see \cite[Sect. 4]{HNspechtfilt},  where the Young module $Y^\lambda$ is constructed. The construction there actually builds $Y^\lambda$ as the left ${\mathcal Y}(\Theta)$-approximation of $S^\lambda$, and can be applied for any $M \in \mo \ks$ in a similar fashion.

By definition the $\Ext$-injectives in the category $\F(\Theta)$ are modules $\F(\Theta) \cap \Y(\Theta)$. We next show the indecomposable objects in this category are exactly the Young modules:

\begin{prop}
\label{extinjectives} The collection of indecomposable modules in $\F(\Theta) \cap \Y(\Theta)$  is exactly $\{Y^\lambda
\mid \lambda \in \Lambda^+(d)\}$.
\end{prop}

\begin{proof} Let $M \in \F(\Theta) \cap \Y(\Theta)$ be indecomposable, so $M$ has a
Specht filtration and $\forall \lambda \in \Lambda^+(d)$:

\begin{eqnarray*}0&=& \Ext^1_{\ks}(S^\lambda,M)\\
 &\cong &\Ext^1_{\ks}(S_{\lambda'},M \otimes \sgn)\\
&\cong& \Ext^1_{S}(\Delta(\lambda'), \G(M \otimes
\sgn)) \mbox{\rm\quad by Prop. \ref{R1Gvanishesondualspechtmodules}(ii)}.\\\end{eqnarray*} 
Thus $\G(M \otimes \sgn)$ is indecomposable and has a good filtration by Prop. \ref{goodweylfiltrationcriterion}. But $M\otimes \sgn$ has a dual Specht filtration so $\G(M \otimes \sgn)$ has a
Weyl filtration by Prop. \ref{propertiesofG}. So  $\G(M \otimes \sgn)$ has both good and Weyl filtration and hence is isomorphic to a tilting module. So $M \otimes \sgn \cong \F(T(\lambda))$ for some $\lambda$, so $M \cong Y^{\lambda '}$ by (\ref{imageofSchur}).
\end{proof}

A similar argument shows the indecomposable Ext-projectives in the category
$\F(\Theta)$ are the twisted Young modules $\{Y^\lambda \otimes
\sgn\}$. We remark that the set $\{(S^\lambda, Y^\lambda)\}$ has
recently been shown to be a {\it stratifying system}, see
\cite{Erdmannstratifsystemsymmetric}. This result together with
Theorem 2.4 of \cite{MarcosMendozaSaenzstratsystviarelative} could
also be used to prove Prop. \ref{extinjectives}.

It would be nice to understand exactly which modules with dual
Specht filtrations are in $\Y(\Theta)$. Of course the simplest module with a dual Specht filtration is just a dual Specht module, and in this case we can give a necessary and sufficient condition for $S_\mu \in \Y(\Theta)$:

\begin{prop}
\label{whenisdualspechtinY} $S_\mu \in \Y(\Theta)$ if and only if
\begin{itemize}
\item[(a)]$\G(S^{\mu'}) \cong \nabla(\mu')$ and

\item[(b)]$R^1\G(S^{\mu'})=0$.
\end{itemize}
\end{prop}
\begin{proof}
Suppose $S_\mu \in \Y(\Theta)$, so
$0=\Ext^1_{\ks}(S^{\tau'},S_\mu)=\Ext^1_{\ks}(S_\tau,S^{\mu'})$
for all $\tau\in \Lambda^+(d)$. From (\ref{fiveterm}) we get that
$\Ext^1(\Delta(\tau),\G(S^{\mu'}))=0 \,\,\,\forall\,\tau$. Thus
$\G(S^{\mu'})$ has a good filtration and maps to $S^{\mu'}$ under
$\F$. Since the multiplicities in a Specht filtration are well-defined we must have $\G(S^{\mu'}) \cong \nabla(\mu')$. Furthermore
$V^{\otimes d}$ has a dual Specht filtration, so $\Ext^1_{\ks}(V^{\otimes d},S^{\mu'})=R^1\G(S^{\mu'})=0$.

Conversely assume $\G(S^{\mu'}) \cong \nabla(\mu')$ and
$R^1\G(S^{\mu'})=0$. Plugging into \eqref{fiveterm} immediately
implies $\Ext^1_{\ks}(S_\tau,S^{\mu'})=0 \, \forall \tau\in \Lambda^+(d)$,
so $S_\mu \in \Y(\Theta)$.
\end{proof}
 It seems a difficult question to determine for which $\mu$ will the conditions of Prop. \ref{whenisdualspechtinY} hold.

Since $\soc (\nabla(\mu')) \cong L(\mu')$,  Props. \ref{propertiesofG}(iv),(v) imply the conditions of
Prop. \ref{whenisdualspechtinY} hold only if $\mu'$ is
$p$-restricted, i.e. if $\mu$ is $p$-regular. A similar argument
shows that $S_\mu$ is detected by Theorem \ref{mainfilttheorem}(i)
only if $\mu$ is $p$-restricted. Thus we have obtained some new information about extensions between Specht and dual Specht modules:

\begin{cor}
\label{corextensionsspechtdualspecht}
\begin{itemize}
\item[] \item[(i)]If $\mu$ is not $p$-regular  then
$\Ext^1_{\ks}(S^\lambda,S_\mu) \neq 0$ for some $\lambda \in \Lambda^+(d)$. \item[(ii)]If $\mu$ is not $p$-restricted then
$\Ext^1_{\ks}(S_\mu,S^\lambda) \neq 0$ for some $\lambda \in \Lambda^+(d)$.

\item[(iii)]If $\mu$ is neither $p$-regular nor $p$-restricted then the
obvious dual Specht filtration of $S_\mu$ is not detected by Theorem \ref{mainfilttheorem}(i) or (ii).
\end{itemize}
\end{cor}
\noindent Setting $\mu=(d)$ in Cor. \ref{corextensionsspechtdualspecht}(ii) gives new information about Specht module cohomology:

\begin{cor}
\label{spechtcohomology}
Let $d \geq p$. Then there exists some $\lambda \in \Lambda^+(d)$ with $H^1(\Sigma_d,S^\lambda) \neq 0$.
\end{cor}
\noindent The previous corollary is in stark contrast to the situation for dual Specht modules where in $\cha p>3$, it is known \cite{BKMdualspecht} that $H^1(\Sigma_d,S_\lambda)=0$ for all $\lambda$ .

\section{A revealing example}
\label{sectionsomeunforunateexamples}
The modules in $\F(\Delta) \cap \F(\nabla)$ are called tilting modules. For each dominant weight $\lambda$ there is a unique indecomposable tilting module $T(\lambda)$ with highest weight $\lambda$. These modules are self-dual and an arbitrary tilting module is a direct sum of them. Since $\F(\nabla)$ and thus $\F(\Delta)$ are closed under tensor products \cite[II.4.21]{jantzenbook2nded}, then so is the collection of tilting modules (as $G$ modules not $S$ modules). Tilting modules also have nice cohomological properties. For example $\Ext^i_S(T(\lambda),T(\mu))=0 \,\,\forall i>0$.

In this section we give an example which shows how these
properties fail for symmetric group modules with both Specht and dual Specht filtrations.
Later, we conjecture that a nice theory may be
salvaged if one assumes additionally that the modules are
self-dual.

Let $p=d=5$. The module category for $k\Sigma_p$ in characteristic
$p$ is completely understood, indeed has only finitely many
indecomposable modules.

\begin{example}
\label{U}Let $U \cong \Omega^2(D_{21^3})$ where $\Omega$ is the
Heller translate. Then $U$ has the following structure:

\setlength{\unitlength}{0.4mm}
\begin{picture}(100,50)(-120,-10)
\put(0,20){$D_{41}$} \put(13,17){\line(1,-1){9}}
\put(53,17){\line(1,-1){9}} \put(43,17){\line(-1,-1){9}}
 \put(40,20){$D_{21^3}$}

\put(20,0){$D_{31^2}$}

\put(60,0){$D_{1^5}$}

\put(-25,10){$U\cong$}
\end{picture}

\end{example}
\noindent The Specht modules in the principal block have the following Loewy structures:
\begin{equation}
\label{Spechtmodulessigma5}
S^5 \cong D_{41},\,\,\, S^{41} \cong \begin{array}{l}
  D_{31^2} \\
  D_{41} \\
\end{array},\,\,\,
S^{31^2} \cong \begin{array}{l}
  D_{21^3} \\
  D_{31^2} \\
\end{array},\,\,\,
S^{21^3} \cong \begin{array}{l}
  D_{1^5} \\
  D_{21^3} \\
\end{array},\,\,\,
S^{1^5} \cong D_{1^5},
\end{equation}
so $U$ has a Specht filtration with
subquotients $S^{1^5}, S^{31^2},$ and $S^5$ and a dual Specht
filtration with subquotients $S_{41}$ and $S_{21^3}$.

\begin{prop}
\label{propertiesofU} Let $U\cong \Omega^2(D_{21^3})$. Then:
\begin{itemize}
\item[(i)]$U \otimes U$ has neither a Specht nor dual Specht
filtration. \item[(ii)]$\Ext^1_{\Sigma_5}(U,U^*) \cong k.$
\item[(iii)]$U$ does not lift to characteristic zero.
\item[(iv)] $U$ has vertex $P \in \operatorname{Syl}_5(\Sigma_5)$. It has source isomorphic to the unique indecomposable $kP$ module of Loewy length three, i.e. $\rad^2(kP)$.
\end{itemize}
\end{prop}
\begin{proof}
We have:
\begin{eqnarray*}
  U \otimes U &\cong& \Omega^2(D_{21^3}) \otimes \Omega^2(D_{21^3}) \\
   &\cong& \Omega^4(D_{21^3} \otimes D_{21^3}) \oplus P \\
   &\cong& \Omega^4(D_{21^3} \oplus D_{41} \oplus D_{32}) \oplus P \\
   &\cong& D_{31^2}\oplus D_{1^5} \oplus P'
\end{eqnarray*}
where $P,P'$ are projective modules. We  used the fact that
$D_{32}$ is projective and that the projective resolutions of the
simple modules are easy to write down. Also the structure $D_{21^3} \otimes
D_{21^3}$ can be easily computed by hand or computer, as $D_{21^3}$ is only three dimensional, the details are left to
the reader. Since projectives have both filtrations and
$D_{31^2}\oplus D_{1^5}$ has neither, we conclude that $U\otimes U$ has neither.

To prove (ii) observe that $U^* \cong \Omega^2(D_{31^2})$ so
$\Ext^1(U,U^*) \cong \Ext^1(D_{21^3},D_{31^2}) \cong k$.
Next notice that if $U$ could be lifted to characteristic zero then
so could $D_{21^3}$, which from \eqref{Spechtmodulessigma5} is not the case.
The last part can also be computed by hand. In particular:
$$\ind_P^{\Sigma_5}\rad^2(kP)
 \cong U \oplus U^* \oplus D_{21^3} \oplus D_{31^2}$$
 
\end{proof}
This example easily generalizes to $\Sigma_p$ for any $p \geq 5$. None of the similar ``pathological" examples which we have constructed have been indecomposable and self-dual. Thus in the next section we consider these modules.

\section{``Tilting" modules for symmetric groups?}
Despite the example in the previous section, we believe that assuming an indecomposable module with both Specht and dual Specht filtrations is also self-dual might allow a complete classification, and we further believe that the geometric tools available in the algebraic group theory could be useful. To further this goal,  we prove that these modules correspond bijectively with a nice class of $G$-modules. Roughly speaking these modules behave like tilting modules under the tilting functor in a way we make precise below.
If our conjecture in the next section is correct, then these modules will be exactly indecomposable self-dual trivial source modules.

\subsection{Ringel duals and tilting functors}
Let $T= \oplus_{\alpha \in \Lambda(n,d)}(\Lambda^\alpha(V))$. Then $T \cong T^\tau$ is a full tilting module and the {\em Ringel dual} of $S(n,d)$ is defined as:

$$S'(n,d) \cong \End_{S(n,d)}(T).$$ The {\em tilting functor} ${\T}: \mo S(n,d) \rightarrow \mo S'(n,d)$ is given by:
$$\T(U)= \Hom_S(T,U).$$ This setup is more thoroughly described in \cite{Donkingtiltingalgebraic}. When $n \geq d$ then Donkin proved \cite[sect. 3]{Donkingtiltingalgebraic} that 

\begin{align}
\label{equationschur=ringeldual}
S(n,d) \cong \End_{k\Sigma_d}(\bigoplus_{\alpha \in \Lambda(n,d)} M^\alpha), &\qquad S'(n,d) \cong \End_{k\Sigma_d}(\bigoplus_{\alpha \in \Lambda(n,d)}M^\alpha \otimes \sgn).
\end{align}

In particular $S(n,d) \cong S'(n,d)$. We let $\tT: \mo S(n,d) \rightarrow \mo S(n,d)$ denote the composition of $\T$ with the functor $\mo S'(n,d) \rightarrow \mo S(n,d)$ arising from the isomorphism. This is the functor  denoted $\tilde{\F}$ in \cite{Donkingtiltingalgebraic}. It is well-known that $\T$, and thus $\tT$, takes modules with good filtrations to modules with Weyl filtrations and interchanges projective and tilting modules. Specifically:

\begin{lemma}\cite[Section 3]{Donkingtiltingalgebraic}
\label{lemmatT}
\begin{itemize}
\item[]
\item[(i)]$\tT(\nabla(\lambda)) \cong \Delta(\lambda')$, $\tT(\F(\nabla)) \subseteq\F(\Delta)$
\item[(ii)]$\tT(T(\lambda)) \cong P(\lambda')$
\item[(iii)]$\tT(P(\lambda)) \cong T(\lambda')$
\item[(iv)]$\tT(I(\lambda))\cong T(\lambda')$
\end{itemize}
\end{lemma}
\noindent In order to relate $\tT$ to modules for $k\Sigma_d$ we need to express $\tT$ in terms of the $\F$ and $\G$:

\begin{prop}
\label{propdonkinstiltingisG(emtimessign)}
Let $U \in \mo S(n,d)$. Then $\tT(U) \cong {\mathcal G}(eU \otimes \sgn)$.
\end{prop}
\begin{proof}
\begin{align*}
\T(U) &\cong \Hom_{S(n,d)}(T,U)\\
&\cong\Hom_{S(n,d)}(U^\tau, T)\\
&\cong \Hom_{k\Sigma_d}((eU)^*, \oplus(M^\alpha \otimes \sgn))\quad \text{since $\G(M^\alpha \otimes \sgn) \cong \Lambda^\alpha(V)$}\\
&\cong \Hom_{k\Sigma_d}(\oplus(M^\alpha \otimes \sgn), eU).
\end{align*}
Now $\T(U)$ is a module for $S(n,d)'$ via its right action on $\oplus(M^\alpha \otimes \sgn)$. To determine $\tT(U)$ we need to compose with the isomorphism from \eqref{equationschur=ringeldual}, which simply takes a map $\phi$ to $\phi \otimes {\rm id}$. Thus we obtain:
\begin{align*}
\tT(U) & \cong \Hom_{k\Sigma_d}(\oplus M^\alpha, eU \otimes \sgn)\\
&\cong \Hom_{k\Sigma_d}(V^{\otimes d}, eU \otimes \sgn)\\
& \cong \G(eU \otimes \sgn)
\end{align*} as desired.
\end{proof}
This alternate description of $\tT$ is all we need to prove:

\begin{theorem}
\label{theoremgivingtiltingcondition}
Let $M \in \mo k\Sigma_d$ be indecomposable and let $U=\G(M)$.
\begin{itemize}
\item[(i)] $M$ has both a Specht and dual Specht filtration if and only if both $U$ and $\tT(U)$ have Weyl filtrations.

\item[(ii)]$M$ is self dual with a Specht (and hence dual Specht) filtration if and only if $U$ has a Weyl filtration and $\tT(U) \cong \tT(U^\tau)$.

\end{itemize}
\end{theorem}

\begin{proof} Notice that $U \in \F(\Delta)$ implies $U^\tau \in \F(\nabla)$ which implies $\T(U^\tau) \in \F(\Delta)$ by Lemma \ref{lemmatT}. Thus the condition $\tT(U) \cong \tT(U^\tau)$ from part (ii) implies (as expected) the condition from part (i) that $\tT(U)$ have a Weyl filtration.

We use Prop. \ref{propdonkinstiltingisG(emtimessign)} repeatedly. First suppose $U$ and $\tT(U)$ both have Weyl filtrations. Then $M \cong eU$ has a dual Specht filtration as does $e\tT(U) \cong M \otimes \sgn$. Thus by \eqref{tensorsignwithspecht}, $M$ has both filtrations. Conversely if $M \in \mo k\Sigma_d$ is indecomposable with both filtrations then so is $M \otimes \sgn$ and setting $U=\G(M)$ we immediately get that both $U$ and $\tT(U)$ have Weyl filtrations by Prop. \ref{propertiesofG}(iv).

Now further suppose that $M \cong M^*$ is indecomposable with both filtrations and let $U=\G(M)$. Then $\tT(U) \cong \G (M \otimes \sgn)$. Similarly $\tT(U^\tau) \cong \G(eU^\tau \otimes \sgn) \cong \G(M^*\otimes \sgn)$ so $\tT(U) \cong \tT(U^\tau)$. The same argument shows that $\tT(U) \cong \tT(U^\tau)$ implies $M$ is self-dual.
\end{proof}

\begin{remark}
\label{remark}
Observe from Lemma \ref{lemmatT} that the tilting modules $U=T(\lambda)$ and the projective modules $U=P(\lambda)$ both have Weyl filtrations and satisfy $\tT(U) \cong \tT(U^\tau)$, so the modules described in Theorem \ref{theoremgivingtiltingcondition} are a simultaneous generalization of tilting modules and projective modules.
\end{remark}

When $n<d$ then $S(n,d)$ and $S(n,d)'$ are usually not isomorphic and so the functor $\tT$ does not exist, but the tilting functor $\T$ does. So one can ask more generally:
\begin{problem}Which $S(n,d)$ modules $U$ (or which $U\in \F(\Delta)$) have the property that $\T(U) \cong \T(U^\tau)$? 
\end{problem}

In \cite{HemmerirreducibleSpecht} we conjectured that indecomposable self-dual $k\Sigma_d$-modules with Specht filtrations would be exactly signed Young modules. Equivalently this conjecture says the modules in Theorem \ref{theoremgivingtiltingcondition} would be the listing modules recently defined by Donkin. A counterexample has recently been found Paget and Wildon:

\begin{example}
\label{examplePaget}
\cite{Pagetspechtdualspecht},\cite{WildonThesis} Let $H \cong C_2^{p} \rtimes \Sigma_p \leq \Sigma_{2p}$ be the normalizer of a fixed-point-free involution in $\Sigma_{2p}$ and let $M \cong \ind_H^{\Sigma_{2p}}k$. Then each indecomposable summand of $M$ is self-dual with a Specht filtration. One of these is not a signed Young module
\end{example}

Specifically Paget proved in \cite{Pagetspechtdualspecht} that the indecomposable summands of $M$ have Specht filtrations, but Wildon had already verified \cite{WildonThesis} that each summand of $M$ is self-dual. It follows by general theory that one of the summands must have vertex a Sylow subgroup $P \leq H$. However $P$ is not conjugate to the Sylow subgroup of any Young subgroup of $\Sigma_{2p}$, so Paget concluded this summand is not a signed Young module. We will discuss this situation further in Section \ref{sectiontrivialsource}.

\section{Irreducible Specht modules}
\label{sectionirreduciblespechtmodules}
Since the irreducible modules for $\ks$ are self-dual, the
irreducible Specht modules give an immediate source of modules
with both Specht and dual Specht filtrations. For $\lambda$
$p$-regular or $p$-restricted the criterion for $S^\lambda$ to be
irreducible has long been known. Only recently Fayers
\cite{FayersirredSpechttypeA} has verified a conjecture of James
and Mathas which settles the problem completely. In this section
we determine which of these modules are detected by the criteria
from Section \ref{filtrationcriterionsection}. 

\begin{lemma}
\label{irreduspechtpregular} Suppose $\mu$ is $p$-regular and
$S^\mu$ is irreducible. Then:
\begin{itemize}
\item[(i)] $S^\mu \cong Y^\mu$.
\item[(ii)]$\Ext_{\ks}(S^\lambda,S_\mu)=0 \, \forall \, \lambda \in \Lambda^+(d)$.
\item[(iii)]$\Ext_{\ks}(S_\mu, S_\lambda)=0 \, \forall \, \lambda \in \Lambda^+(d)$.
\end{itemize}
\end{lemma}
\begin{proof}
Since $\mu$ is $p$-regular then $S^\mu \cong D^{\mu}$. Also $[Y^\lambda:D^\lambda]=1$. But $Y^\lambda$ is self-dual and
indecomposable so $Y^\lambda \cong D^\lambda \cong S^\lambda$.
Parts (ii) and (iii) follow from  the fact \cite[4.1.1]{HNspechtfilt} that $\Ext^1(S^\lambda,
Y^\mu)=0$ for $p>3$ and $S^\mu$ is self-dual.
\end{proof}
Notice parts (ii) and (iii) are saying that both the Specht and
dual Specht filtrations of $S^\lambda$ are detected by Theorem
\ref{mainfilttheorem}. Similarly we have:

\begin{lemma}
\label{irreduspechtprestr} Suppose $\mu$ is $p$-restricted and
$S^\mu$ is irreducible. Then:
\begin{itemize}
\item[(i)] $S^\mu \cong Y^{\mu'} \otimes \sgn$.
\item[(ii)]$\Ext_{\ks}(S_\lambda,S_\mu)=0 \, \forall \, \lambda
\in \Lambda^+(d)$. 
\item[(ii)]$\Ext_{\ks}(S_\mu, S^\lambda)=0 \,
\forall \, \lambda \in \Lambda^+(d)$.
\end{itemize}
\end{lemma}
\begin{proof} This follows from Lemma \ref{irreduspechtpregular}
immediately  since $S^\mu \otimes \sgn \cong S_{\mu'}\cong
S^{\mu'}$ and $\mu'$ is $p$-regular.
\end{proof}

In stark contrast to the previous two lemmas the next result shows
that irreducible Specht modules which are neither $p$-restricted
nor $p$-regular are {\it never} detected by Theorem
\ref{mainfilttheorem}:

\begin{prop}
\label{testsmissnonpregpsignirrSpec}
 Suppose $\mu$ is neither
$p$-regular nor $p$-restricted and suppose $S^\mu$ is
irreducible. Then $S^\mu$ does not satisfy any of the criteria
in Theorem \ref{mainfilttheorem}.
\end{prop}
\begin{proof}Suppose $S^\mu$ satisfies the criterion in Theorem \ref{mainfilttheorem}(ii), so $\Ext^1_{\ks}(S^\lambda,S^\mu)=0 \,\, \forall
\lambda.$ Then $S^\mu \cong Y^\mu$ by Prop.
\ref{extinjectives}. So $\G(S^\mu) \cong \Delta(\mu) \cong P(\mu)$.
Thus $\Delta(\mu)\cong L(\mu)$ is an irreducible projective module
and $\mu$ must be $p$-restricted, contradicting our assumption. Thus $S^\mu$ cannot satisfy
criterion (ii). Similar arguments handle the other three criteria.
\end{proof}

\section{Trivial Source modules} 
\label{sectiontrivialsource}
Recall that the Young modules $Y^\lambda$ are precisely the
indecomposable summands of the permutation modules $M^\lambda
\cong {\rm Ind}_{\Sigma_\lambda}^{\Sigma_d}k$. For $\lambda \vdash
a$ and $\mu \vdash b$ with $a+b=d$ define the {\it signed
permutation module}:
$$M^{(\lambda \mid \mu)} \cong \ind_{\Sigma_\lambda \times
\Sigma_\mu}^{\Sigma_d}k \boxtimes \sgn.$$ Indecomposable summands of signed permutation modules are called {\em signed Young modules}. Notice the class of
signed Young modules includes all ordinary Young modules $Y^\lambda$ and
twisted Young modules $Y^\lambda \otimes \sgn$ but, when $d \geq 2p$,  includes other modules
as well.
Signed Young modules give a large class of self-dual modules with
both Specht and dual Specht filtrations. We once believed this was exactly the class of indecomposable self-dual modules with Specht filtrations. However Example \ref{examplePaget} shows there are more.

The signed Young modules have trivial source and vertex isomorphic to a Sylow subgroup of some Young subgroup of $\Sigma_d$. However, in general there are many conjugacy classes of $p$-subgroups of $\Sigma_d$ which are not of this form, and many more trivial source modules, for example the module $M$ in Example \ref{examplePaget}. On the other hand the module $U$ from Prop. \ref{propertiesofU} has both Specht and dual Specht filtrations, and is not a trivial-source module. But it is not self-dual either!  Undaunted, we dare to hazard another conjecture:

\begin{conj}
\label{wildconjecture}
Indecomposable, self-dual $k\Sigma_d$ modules with Specht filtrations are trivial source modules.
\end{conj}
If this conjecture were true it would have the nice consequence that there are only finitely many indecomposable self-dual $k\Sigma_d$ modules with Specht filtrations.

As evidence for this we showed:

\begin{theorem}\cite[Theorem 4.2]{HemmerirreducibleSpecht}
\label{irrspechtmodulesaresignedyoung} Suppose $S^\mu$ be
irreducible. Then $S^\mu$ is a signed Young module, and hence has trivial source.
\end{theorem}

\section{Final remarks and problems} Even though many of the results and the machinery of this paper depend on the characteristic being at least five, it is possible that Conjecture \ref{wildconjecture} holds in characteristics two and three. For example in characteristic two, $d=4$, there are two irreducible $\Sigma_4$-modules and both of them are isomorphic to Specht modules and both have trivial source. Thus \emph{every} $\Sigma_4$ module has both a Specht and dual Specht filtration. For Conjecture \ref{wildconjecture} to hold here, it would require $k\Sigma_4$ to have only finitely many indecomposable self-dual modules, even though the algebra has wild type!

We close with a list of problems for further study.

\begin{problem} Suppose $M$ and $N$ are indecomposable self-dual modules with Specht filtrations. Is $M \otimes N$ a direct sum of such modules?
\end{problem}

We have been unable to answer this even for $M \otimes M$ where $M$ is from Example \ref{examplePaget}.

\begin{problem} Must indecomposable, trivial source modules have either Specht or dual Specht filtrations? More generally, must $k\Sigma_d$ modules which are $p$-modular reductions of some indecomposable ${\mathcal O}\Sigma_d$-lattice have Specht or dual Specht filtrations?
\end{problem}
For $\Sigma_p$ this is known, the trivial source modules are exactly the Young  and twisted Young modules, together with the Specht modules $S^{(p-a,1^a)}$ and their duals $S_{(p-a,1^a)}$ for $a$ even \cite[4.3]{ErdmannYoungsymmetricgroups}.

\begin{problem}\label{remarkextvanishishingsignedYoung} Suppose $M$ and $N$ are indecomposable self-dual modules with Specht filtrations and suppose $p>3$.
Is $\Ext^i_{\ks}(M,N)=0$ for $ 1 \leq i \leq p-3$? 
\end{problem}

\noindent This is known to be the case for signed Young modules \cite[Lemma 2.2(iii)]{HemmerirreducibleSpecht}.

Notice that $U$ in Example \ref{U} does not lift to characteristic zero, by Prop. \ref{propertiesofU}(iii). This is true of other indecomposable modules we have found which have both Specht and dual Specht filtrations but are not self dual. The signed Young modules do, however, lift to characteristic zero, as does any trivial source module.

\begin{problem} Suppose $M$ and $N$ are indecomposable self-dual modules with Specht filtrations and suppose $p>3$. Do $M$ and $N$ lift to characteristic zero?
\end{problem}
 If, for example, Problem \ref{remarkextvanishishingsignedYoung} were answered affirmatively, it would in particular show that $\Ext^2_{\ks}(M,M)=0$. This $\Ext^2$-vanishing would \cite[3.7.7]{Bensonreptheoryvolume1} imply that $M$ lifts to characteristic zero.  We remark that the module itself lifting to an integral representation is a weaker condition than the entire filtration lifting.

Notice that $U$ from Example \ref{U} does not self-extend, i.e. $\Ext^1_{\ks}(U,U)=0$. We know of no examples of indecomposable modules with Specht and dual Specht filtrations that self-extend. 

\begin{problem}Suppose $U \in \mo \ks$ is indecomposable with both Specht and dual Specht filtrations. Must $\Ext^1_{\ks}(U,U)=0$?
\end{problem}

\bibliographystyle{amsplain}
\bibliography{references}
\end{document}